\newcount\notenumber

\def\note{\advance\notenumber by 1
\footnote{$^{(\the\notenumber)}$}}

\def\C{{\bf C}}
   
\def\Gal{{\rm Gal}}

\def\OS{{\cal O}_S}

\def\Q{{\bf Q}}
\def\R{{\bf R}}
\def\N{{\bf N}}
\def\Z{{\bf Z}}

\def\x{{\bf x}}

\hsize = 15truecm
\vsize = 22truecm

\hoffset = 0.6truecm
\voffset = 0.7truecm


\font\title=cmr10 scaled 1200


\centerline{\title On the rational approximations to the powers of an
algebraic number}

\centerline{\it Solution of two problems of Mahler and Mend\`es France}\bigskip

\centerline{Pietro Corvaja\qquad Umberto Zannier}\bigskip

\noindent{\bf Abstract.} About fifty years ago Mahler [M] proved that {\it if
$\alpha>1$ is   rational but not an integer and if $0<l<1$, then the fractional
part of $\alpha^n$ is $>l^n$ apart from a finite set of integers $n$ depending
on $\alpha$ and $l$}. Answering completely a question of Mahler, we show (Thm. 1) that 
the same conclusion holds for all algebraic numbers which are not $d$-th roots of
  Pisot numbers. By related methods we also answer a question of 
Mend\`es France, characterizing completely the quadratic irrationals $\alpha$ such that the
continued fraction for $\alpha^n$ has period length tending to infinity (Thm. 2).\bigskip

\noindent{\bf \S 1 Introduction.} About fifty years ago Mahler [M] proved that {\it if
$\alpha>1$ is   rational but not an integer and if $0<l<1$, then the fractional
part of $\alpha^n$ is $>l^n$ apart from a finite set of integers $n$ depending
on $\alpha$ and $l$}.
His proof used a $p$-adic version of Roth's Theorem, as in previous work by Mahler and
especially by Ridout.\smallskip

At the end of that paper Mahler pointed out that the conclusion does not hold
if $\alpha$ is a suitable algebraic number, as e.g. ${1\over 2}(1+\sqrt 5)$;
of course a counter-example is provided by any {\it Pisot} number, i.e. a real
algebraic integer $\alpha>1$ all of whose conjugates different
from $\alpha$ have absolute value $<1$ (note that rational integers $>1$ 
are Pisot numbers according to our definition).
Mahler also added that ``It would be of some interest to know which algebraic
numbers have the same property as..." the rationals in his theorem.

Now, it seems that even replacing Ridout Theorem with the modern versions of
Roth's Theorem, valid for several valuations and approximations in any given
number field, the method of Mahler does not lead to a complete solution to his
question.

One of  the objects of the present paper is to answer completely Mahler's question; our
methods will involve a suitable version of the Schmidt Subspace Theorem,
which may be considered as a multi-dimensional extension of the
mentioned results by Roth, Mahler and Ridout. We state at once our first
theorem, where as usual we denote by $||x||$ the distance of the real
number $x$ from the nearest integer.\medskip

\noindent {\bf Theorem 1.}\ {\it Let $\alpha >1$ be a real algebraic number 
and let $0<l<1$. Suppose that $||\alpha^n||<l^n$ for infinitely many
 natural numbers $n$. Then some power $\alpha^n$  is  a Pisot number. 
In particular, $\alpha$ is an algebraic integer}.\medskip

We remark that the conclusion is   best-possible, since, if it is
true, then conversely we have   $\|\alpha^n\|\ll l^n$ 
for all $n$ in a suitable arithmetic progression,
  where $l$ is the maximum absolute value
of the conjugates of $\alpha$ different from $\alpha$.  
Here, Mahler's example with the Golden ratio is  
typical. 

Also, the conclusion is not generally true without the assumption that $\alpha$ is algebraic;
see for this the Appendix.\smallskip

The  present application of the Subspace Theorem seems different from previous ones,
and occurs in Lemma 3 below. Related methods actually enable us to answer as well a question
raised by Mend\`es France about the length of the periods of the continued fractions for
$\alpha^n$, where $\alpha$ is now a quadratic irrational; this   appears as  Problem 6 in
[MF]. We shall prove, more generally: \smallskip

\noindent {\bf Theorem 2}. {\it Let $\alpha>0$ be a real quadratic irrational.
If $\alpha$ is neither the square root of a rational number, nor
a unit in the ring of integers of $\Q(\alpha)$ then 
the period of the continued fraction   for $\alpha^n$   tends  
to infinity with $n$. If $\alpha$ is the square root of a rational number 
the period length of the continued fraction for
$\alpha^{2n+1}$ tends to infinity. If $\alpha$ is a unit, the period
of the continued fraction for $\alpha^n$ is bounded.
}
\medskip

Clearly, if $\alpha$ is the square root of a rational number,
then the continued fraction for  $\alpha^{2n}$ is finite\note{its length tends to infinity by
a result of Pourchet, proved in greater generality in [CZ]}, so Theorem 2 gives a complete
answer to the problem of Mend\`es France.

The main tool in the proof of both theorems is the following new lower
bound for the fractional parts of $S$-units in algebraic number fields.
We first need a definition:
\medskip

\noindent{\bf Definition}. We call a (real) algebraic number $\alpha$   a 
{\it pseudo-Pisot} number if:

\item{$(i)$} $|\alpha|>1$ and all its conjugates have (complex)
absolute value strictly less then $1$;

\item{$(ii)$} $\alpha$ has an integral trace: ${\rm Tr}_{\Q(\alpha)/\Q}(\alpha)\in\Z$.

Of course, pseudo-Pisot numbers are ``well approximated" by their trace,
hence are good candidates for having a small fractional part 
compared to their height. The algebraic integers among the pseudo-Pisot 
numbers are just the usual Pisot numbers. We shall prove
\medskip

\noindent{\bf Main Theorem}. {\it Let $\Gamma\subset\overline{\Q}^\times$ be a 
finitely generated multiplicative group of algebraic numbers, let 
$\delta\in\overline{\Q}^\times$ be a non zero algebraic number
 and let $\epsilon>0$ be fixed.
Then there are only finitely many pairs $(q,u)\in\Z\times \Gamma$
with $d=[\Q(u):\Q]$
such that $|\delta qu|>1$, $\delta qu$ is not a pseudo-Pisot number and
$$
0<\|\delta qu\|<H(u)^{-\epsilon}q^{-d-\epsilon}.\eqno(1.1)
$$
}
\medskip

Note again  that, conversely, starting with a Pisot number $\alpha$ and taking 
$q=1,u=\alpha^n$ for $n=1,2,\ldots$ produces an infinite 
sequence of solutions to $0<\|qu\|<H(u)^{-\epsilon}$ for a 
suitable $\epsilon>0$.

\bigskip

The above Main Theorem can be viewed as a Thue-Roth inequality with
``moving target", as the Theorem in [CZ], where we considered
quotients of power sums with integral roots instead of elements of
a finitely generated multiplicative group. The main application of
the Theorem in [CZ] also concerned continued fractions, as for our Theorem 2.


\bigskip

\noindent{\bf \S 2. Proofs}. We shall use the following notations: let $K$ be
a number field, embedded in $\C$, Galois over $\Q$.
We denote by $M_K$ (resp. $M_\infty$) be the set of places 
(resp. archimedean places) of $K$;
for each place $v$ we denote by $|\cdot|_v$  the absolute value
corresponding to  $v$, normalized {\it with respect to} $K$;
by this we mean that if $v\in M_\infty$ then there exists an
automorphism  $\sigma\in\Gal(K/\Q)$  of $K$ such that, for all $x\in K$
$$
|x|_v=|\sigma(x)|^{d(\sigma)\over [K:\Q]}
$$
where $d(\sigma)=1$ if $\sigma(K)=K\subset\R$ and $2$ otherwise
(note that $d(\sigma)$ is now constant since $K/\Q$ is Galois).
Non-archimedean absolute values are normalized accordingly, 
so that the product formula holds and the absolute Weil height reads
$$
H(x)=\prod_{v\in M_k}\max\{1,|x|_v\}.
$$
For a vector $\x=(x_1,\ldots,x_n)\in K^n$ and a place $v\in M_K$
we shall denote by $\|\x\|_v$ the $v$-norm of $\x$:
$$
\|\x\|_v:=\max\{|x_1|_v\ldots,|x_n|_v\},
$$
and by $H(\x)$ the projective height
$$
H(\x)=\prod_{v\in M_K}\max\{|x_1|_v,\ldots,|x_n|_v\}.
$$
\smallskip

We begin by proving the Main Theorem. First of all 
notice that, by enlarging if necessary the multiplicative
group $\Gamma$,  we can reduce to the following
situation: $\Gamma\subset K^\times $ is the group of $S$-units:
$$
\Gamma=\OS^\times=\{u\in K\quad : \quad |u|_v=1\quad {\rm for\ all}\ v\not\in S\} 
$$
of a suitable number field $K$, Galois over $\Q$, with respect to
a suitable finite set of places $S$, containing the
archimedean ones and stable under Galois conjugation.

\medskip

Our first lemma can be easily deduced from a theorem
of Evertse, which in turn was obtained as an application of the
already mentioned Subspace Theorem; we prefer to give a  
proof for completeness.
\medskip

\noindent{\bf Lemma 1}. {\it Let $K, S$ be as before, $\sigma_1,\ldots,\sigma_n$
distinct automorphisms of $K$, $\lambda_1,\ldots,\lambda_n$ non zero 
elements of $K$, $\epsilon>0$ be a positive
real number, $w\in S$ be a distinguished place.
Let $\Xi\subset\OS^\times$ be an infinite
set of  solutions $u\in\OS^\times$ of the inequality
$$
|\lambda_1\sigma_1(u)+\ldots+\lambda_n\sigma_n(u)|_w<
\max\{|\sigma_1(u)|_w,\ldots,|\sigma_n(u)|_w\}\cdot H(u)^{-\epsilon}.
$$
Then there exists a non trivial   linear relation  of the form 
$$
a_1 \sigma_1(u)+\ldots+a_n\sigma_n(u)=0, \qquad a_i\in K,
$$
satisfied by infinitely many elements of  $\Xi$.}
\medskip

\noindent{\it Proof}. Let $\Xi$ be as in the Lemma. Going to an infinite subset of $\Xi$ we
may assume that $|\sigma_1(u)|_w=\max\{|\sigma_1(u)|_w,\ldots,|\sigma_n(u)|_w\}$  for all the
involved $u$'s.  
Let us consider, for each $v\in S$, $n$ linear forms $L_{v,1},\ldots,L_{v,n}$
in $n$ variables $\x=(x_1,\ldots,x_n)$ as follows:   put
$$
L_{w,1}(\x):=\lambda_1x_1+\ldots+\lambda_n x_n
$$
while for $(v,i)\in S\times\{1,\ldots,n\}$, with 
$(v,i)\neq (w,1)$, put $L_{v,i}(\x)=x_i$. Note that the linear forms 
$L_{v,1},\ldots,L_{v,n}$ are indeed linearly independent. Now put
$\x=(\sigma_1(u),\ldots,\sigma_d(u))\in{\OS^\times}^d$ and
consider the double product
$$
\prod_{v\in S}\prod_{i=1}^n{|L_{v,i}(\x)|_v\over \|\x\|_v}.
$$
By multiplying and dividing by $|x_1|_w=|\sigma_1(u)|_w$, and using the fact that
the coordinates of $\x$ are $S$-units, we obtain the equality
$$
\prod_{v\in S}\prod_{i=1}^n{|L_{v,i}(\x)|_v\over \|\x\|_v}=
|\lambda_1\sigma_1(u)+\ldots+\lambda_n\sigma_n(u)|_w\cdot|\sigma_1(u)|_w^{-1}
\cdot H(\x)^{-n}.
$$
Since $u\in \Xi$, we have   
$$ 
|\lambda_1\sigma_1(u)+\ldots+\lambda_n\sigma_n(u)|_w\cdot |\sigma_1(u)|_w^{-1}
<H(u)^{-\epsilon},
$$ 
then 
$$
\prod_{v\in S}\prod_{i=1}^n{|L_{v,i}(\x)|_v\over \|\x\|_v}<
H(\x)^{-n}H(u)^{-\epsilon}.
$$
The height of the point $\x=(\sigma_1(u),\ldots,\sigma_d(u))$ is easily 
compared with the height of $u$ by the estimate 
$$
H(\x)\leq H(u)^{[K:\Q]}
$$
hence the above upper bound for the double product also gives
$$
\prod_{v\in S} \prod_{i=1}^n{|L_{v,i}(\x)|_v\over \|\x\|_v}<
H(\x)^{-n-\epsilon/[K:\Q]}.
$$
Now, an application of the Subspace Theorem in the form  given in e.g.
[S, Theorem 1D$'$] gives the desired result.

\medskip

Our next tool is a very special case of
the so called unit-equation-theorem, proved by
Evertse and   van der Poorten-Schlickewei. It also rests on the
Subspace Theorem:
\medskip

\noindent{\bf Lemma 2}. {\it Let $K,S,\sigma_1,\ldots,\sigma_n$ be as before,
$\epsilon>0$  a real number
 $a_1,\ldots,a_n$
be non zero elements of $K$. Suppose that $\Xi\subset \OS^\times$
is an infinite set of solutions to the equation
$$
a_1\sigma_1(u)+\ldots+a_n\sigma_d(u)=0.
$$
Then there exist two distinct indices $i\neq j$, two non zero elements
$a,b\in K^\times$ and an infinite subset
$\bar{\Xi}\subset\Xi$ of $\Xi$ such that for all $u\in\bar{\Xi}$
$$
a\sigma_i(u)+b\sigma_j(u)=0.
$$
}

For a proof of this lemma, see  [S, Chap. 4].
\medskip

Our last lemma is the key of the proof of the Main Theorem; its
proof is once more based on  the Subspace Theorem.
\medskip

\noindent{\bf Lemma 3}. {\it Let $K,S$ be as before, $k\subset K\cap\R$  a 
(real) subfield of $K$, of degree $d$ over $\Q$,
$\delta\in K^\times$ a nonzero element of $K$. Let $\epsilon>0$ 
be given. Suppose we have an infinite sequence $\Sigma$ of
points  $(q,u)\in \Z\times(\OS^\times\cap k)$  
such that $|\delta qu|>1$, $\delta qu$ is not a pseudo-Pisot number and
$$
0<\|\delta q u\|<H(u)^{-\epsilon}q^{-d-\epsilon}.\eqno(2.1)
$$
Then there exists a proper subfield $k^\prime\subset k$, an element 
$\delta^\prime\in k^\times$ and an infinite subsequence $\Sigma^\prime\subset\Sigma$ 
such that for all $(q,u)\in \Sigma^\prime$, $u/\delta^\prime\in k^\prime$.
}
\medskip

Note that Lemma 3 gives a finiteness result in the case
$k=\Q$, since the rational field $\Q$ admits no
proper subfields.
\medskip

\noindent{\it Proof of Lemma 3}. Let us suppose that the hypotheses
of the Lemma are satisfied, so in particular $\Sigma$ is an infinite sequence
of solutions of (2.1). We begin by observing that, by  Roth's theorem
[S, Theorem 2A], in any such   infinite sequence $u$ cannot be fixed; so we  
have $H(u)\rightarrow\infty$ in the set $\Sigma$. 
Let $H:=\Gal(K/k)\subset\Gal(K/\Q)$ be the subgroup   fixing  $k$        
(note that $H$ contains the complex conjugation). Let $\{\sigma_1,\ldots,\sigma_d\}$
($d=[k:\Q]$) be a (complete) set of representatives for the left cosets
of $H$ in  $\Gal(K/\Q)$, containing the identity $\sigma_1$.
Each automorphism $\rho\in \Gal(K/\Q)$ defines an archimedean
 valuation on $K$ by the formula
$$
|x|_{\rho}:= |\rho^{-1}(x)|^{d(\rho)\over [K:\Q]}\eqno(2.2)
$$
where as usual $|\cdot|$ denotes the usual complex absolute value.
Two distinct automorphisms $\rho_1,\rho_2$ define the same valuations if
and only if $\rho_1^{-1}\circ\rho_2$ is the complex conjugation. Note that
in this case $\rho_1,\rho_2$ coincide on $k$.
Let now $(q,u)\in \Sigma$ be a solution of (2.1) and let $p\in\Z$ be the
nearest integer to $\delta qu$. 
Then for each $\rho\in\Gal(K/\Q)$ we have, with the notation of (2.2), 
$$
\|\delta qu\|^{d(\rho)\over [K:\Q]}=
|\delta qu-p|^{d(\rho)\over [K:\Q]}=|\rho(\delta)\rho(qu)-p|_\rho.\eqno(2.3)
$$
Let, for each $v\in M_\infty$, $\rho_v$ be an
automorphism defining the valuation $v$ according to the rule (2.2): $|x|_v:=|x|_{\rho_v}$;
then the set $\{\rho_v\, |\, v\in M_\infty\}$ represents the left cosets
of the subgroup generated by the complex conjugation in $\Gal(K/\Q)$.
Let $S_i$, for $i=1,\ldots,d$, be the subset of $M_\infty$ formed by those
valuations $v$ such that $\rho_v$ coincides with $\sigma_i$ on $k$;
note that $S_1\cup\ldots \cup S_d=M_\infty$.
We take the product of the terms in   (2.3)
for $\rho$ running over the set $\{\rho_v\, |\, v\in M_\infty\}$:
this corresponds to taking the product over all archimedean valuations.
Then we obtain
$$
\prod_{v\in M_\infty} |\rho_v(\delta)\rho_v(qu)-p|_{v}=
\prod_{i=1}^d  \prod_{v\in S_i} |\rho_v(\delta)\sigma_i(qu)-p|_{v}.
$$
By (2.3) and the well-known formula $\sum_{v\in M_\infty} d(\rho_v)=[K:\Q]$
it follows that  
$$
\prod_{i=1}^d  \prod_{v\in S_i} |\rho_v(\delta)\sigma_i(qu)-p|_{v}=
\|\delta qu\|.   \eqno(2.4)
$$
Now, let us define for each $v\in S$ a set of $d+1$ linearly independent linear
forms in $d+1$ variables $(x_0,x_1,\ldots,x_d)$ in the following way:
for an archimedean valuation  $v\in S_i$ ($i=1,\ldots,d$) put
$$
L_{v,0}(x_0,x_1,\ldots,x_d)=x_0-\rho_v(\delta) x_i
$$
and for all $v\in S\setminus M_\infty$ or $0<j\leq d$ put
$$
L_{v,j}(x_0,x_1,\ldots,x_d)=x_j.
$$
Plainly the forms $L_{v,0},\ldots,L_{v,d}$ are independent for each $v\in S$.
Finally, let $\x\in K^{d+1}$ be the point
$$
\x=(p,q\sigma_1(u),\ldots,q\sigma_d(u))\in K^{d+1}.
$$    
Let us estimate the double product
$$
\prod_{v\in S}\prod_{j=0}^d {|L_{v,j}(\x)|_v\over\|\x\|_v}.\eqno(2.5)
$$
Using the fact that $L_{v,j}(\x)=q\cdot\sigma_j(u)$ for $j\geq 1$
and that the $\sigma_j(u)$
are $S$-units, we obtain, from the product formula,
$$
\prod_{v\in S}\prod_{j=1}^d  |L_{v,j}(\x)|_v\leq \prod_{v\in M_\infty}\prod_{j=1}^d |q|_v
= |q|^d.\eqno(2.6)
$$
Since the coordinates of $\x$ are $S$-integers,
the product of the denominators in (2.5) is $H(\x)^{d+1}$; then in view of
(2.4) and (2.6)
$$
\prod_{v\in S}\prod_{j=0}^d {|L_{v,j}(\x)|_v\over\|\x\|_v}\leq
H(\x)^{-d-1}\cdot |q|^d\cdot  \| \delta qu\|\leq
H(\x)^{-d-1} (qH(u))^{-\epsilon},
$$
the last inequality being justified by the fact that $(q,u)\in\Sigma$, 
hence (2.1) holds.
Since $H(\x)\leq |q|\cdot |p|\cdot H(u)^{d}$ and $|p|\leq |\delta qu|+1\leq
|q|\cdot |\delta|\cdot H(u)^d+1\leq |q|H(u)^{2d}$ for all but finitely many 
pairs $(q,u)\in\Sigma$ (recall that $H(u)\rightarrow\infty$ for $(q,u)\in\Sigma$),
 we have $|q|^2 H(u)^{3d}\geq  H(\x)$ so  the last displayed
inequality gives
$$
\prod_{v\in S}\prod_{j=0}^d {|L_{j,v}(\x)|_v\over\|\x\|_v}\leq
H(\x)^{-d-1-(\epsilon/3d)}.
$$
An application of the Subspace Theorem in the form given e.g. in
[S, Theorem 1D$'$] implies the existence of a hyperplane containing
infinitely many points $\x=(p,q\sigma_1(u),\ldots,q\sigma_d(u))$.
We then obtain a non trivial linear relation of the form
$$
a_0p+a_1q\sigma_1(u)+\ldots+a_dq\sigma_d(u)=0, \qquad a_i\in K, \eqno(2.7)
$$
satisfied by infinitely many pairs $(q,u)\in\Sigma$.  Our next goal is to prove the
following:  

{\bf Claim}: {\it there exists one such non-trivial relation with vanishing coefficient $a_0$, 
i.e. one involving only the conjugates of $u$}. 

We prove the Claim. Rewriting the above linear relation, if $a_0\neq 0$,  we obtain
$$
p=-{a_1\over a_0}q\sigma_1(u)-\ldots-{a_d\over a_0}q\sigma_d(u).\eqno(2.8)
$$
Suppose first that for some index 
 $j\in\{2,\ldots,d\}$,  $\sigma_j(a_1/a_0)\neq a_j/a_0$; then by applying
the automorphism $\sigma_j$ to both sides of (2.8)
and subtracting term-by-term from  (2.8) to eliminate $p$, we obtain a linear relation 
involving only the terms $\sigma_1(u),\ldots,\sigma_d(u)$. Such a relation
is non trivial since the coefficient of $\sigma_j(u)$ becomes
$\sigma_j(a_1/a_0)- a_j/a_0$. Hence we have proved our Claim in this case.

Therefore, we may and  shall assume that
$a_j/a_0=\sigma_j (a_1/a_0)$ for all $j$,   so in particular all coefficients 
$a_j/a_0$ are nonzero. Let us then rewrite (2.8) in a simpler form as
$$
p=q(\sigma_1(\lambda)\sigma_1(u)+\ldots+\sigma_d(\lambda)\sigma_d(u)) 
\eqno(2.9)
$$
with $\lambda=-(a_1/a_0)\neq 0$.
Suppose now that $\lambda$ does not belong to $k$. Then there exists an automorphism
$\tau\in H$ with $\tau(\lambda)\neq\lambda$. (Recall that $H$ is the subgroup
of $\Gal(K/\Q)$ fixing $k$ and that $\{\sigma_1=id,\sigma_2,\ldots,\sigma_d\}$
is a complete set of representatives of left cosets of $H$).
By applying the automorphism $\tau$ to both sides of (2.9) 
and subtracting term-by-term from  (2.9) to eliminate $p$, we obtain the linear relation
$$
(\lambda-\tau(\lambda))\sigma_1(u)+
(\sigma_2(\lambda)\sigma_2(u)-\tau\circ\sigma_2(\lambda)\tau\circ\sigma_2(u))+\ldots+
(\sigma_d(\lambda)\sigma_d(u)-\tau\circ\sigma_d(\lambda)\tau\circ\sigma_d(u))=0.
$$
Note that $\tau\circ\sigma_j$  coincides on $k$ with some $\sigma_i$. Note also that since
$\tau\in H$ and $\sigma_2,\ldots,\sigma_d\not\in H$, no   
$\tau\circ\sigma_j$ with $j\geq 2$ can belong to $H$. Hence  the above relation can be again
written as a linear combination of the $\sigma_i(u)$; in such expression  the coefficient of
$\sigma_1(u)$ will remain $\lambda-\tau(\lambda)$ and will therefore be $\neq 0$, so we obtain
a non trivial relation among the $\sigma_i(u)$, as claimed.

Therefore we may and shall suppose that $\lambda\in k$ and write (2.9)
in the simpler form
$$
p=q\cdot  {\rm Tr}_{k/\Q}(\lambda u)={\rm Tr}_{k/\Q}(q\lambda u).\eqno(2.10)
$$
After adding $-\delta qu$ to both sides in (2.9), 
recalling that $(q,u)$ is a solution to
(2.1) and that $\sigma_1$ is the identity, we obtain
$$
|p-\delta qu|=|(\lambda-\delta) q\sigma_1(u)+ q\sigma_2(\lambda)\sigma_2(u)+\ldots
+q\sigma_d(\lambda)\sigma_d(u)|<q^{-d-\epsilon}H(u)^{-\epsilon}\leq 
q^{-1}H(u)^{-\epsilon}. 
$$
In particular,
$$
|(\lambda-\delta)u+\sigma_2(\lambda)\sigma_2(u)+\ldots +\sigma_d(\lambda)\sigma_d(u)|
< q^{-1} H(u)^{-\epsilon}.\eqno(2.11)
$$
We want to apply Lemma 1. We distinguish two cases:

{\it First case}: $\lambda=\delta$ (in particular $\delta\in k$).

In this case the algebraic number $q\delta u=q\lambda u$ has an integral trace.
Since by assumption it is not a pseudo-Pisot number, the maximum modulus
of its conjugates $|\sigma_2(q\lambda u)|,\ldots,|\sigma_d(q\lambda u)|$
is $\geq 1$. This yields
$$
 \max\{|\sigma_2(u)|,\ldots,|\sigma_d(u)|\}
\geq q^{-1}\max\{|\sigma_2(\lambda)|,\ldots,|\sigma_d(\lambda)|\}^{-1}.
$$
Hence from (2.11) we deduce that infinitely many pairs $(q,u)\in\Sigma$
satisfy the inequality of Lemma 1, where   $w$ is the  archimedean place 
associated to the given embedding   $K\hookrightarrow \C$, $n=d-1$,
$\lambda_i=\sigma_{i+1}(\lambda)$ and with $\sigma_2,\ldots,\sigma_d$
instead of $\sigma_1,\ldots,\sigma_n$. The conclusion of Lemma 1 
provides what we claimed.

{\it Second case}. $\lambda\neq \delta$. 

In this case the first term does appear in (2.11). Since we supposed
$|q\delta u|>1$, we have 
$$
\max\{|\sigma_1(u)|,\ldots,|\sigma_d(u)|\}\geq |u|>|\delta|^{-1}\cdot q^{-1}
$$
so we can again apply Lemma 1 (with the same place $w$ as in the first case,
$n=d$, $\lambda_1=(\lambda-\delta),\lambda_2=\sigma_2(\lambda),\ldots,
\lambda_d=\sigma_d(\lambda)$) and conclude as in the first case. 
 
This finishes the proof of the Claim, i.e.
a relation of the kind (2.7) without the term $a_0p$ is satisfied for all
pairs $(q,u)$ in an infinite set $\bar{\Sigma}\subset\Sigma$. \smallskip


We can now apply the unit theorem of Evertse
and van der Poorten-Schlickewei in the form of Lemma 2
 which implies that a non trivial relation of
the form $a\sigma_j(u)+b\sigma_i(u)=0$ for some $i\neq j$ and $a,b\in K^\times$ 
is satisfied for $(q,u)$ in an infinite subset $\Sigma^\prime$ of $\bar{\Sigma}$. 
We rewrite it as
$$
-\sigma_i^{-1}\left({a\over b}\right)\cdot (\sigma_i^{-1}\circ\sigma_j)(u)=u.
$$
Then for any two solutions $(q^{\prime},u^\prime), 
(q^{\prime\prime},u^{\prime\prime})\in{\Sigma^\prime}$, 
the element $u:=u^\prime/u^{\prime\prime} \in k$ is fixed by
the automorphism $\sigma_i^{-1}\circ\sigma_j\not\in H$, thus
$u$ belongs to the proper subfield $k^\prime:=K^{<H,\sigma_i^{-1}\circ\sigma_j>}$
of $k$. In other words, if we let $(q^\prime,\delta^\prime)$ be any solution, then
infinitely many solutions are of the form $(q^{\prime\prime},v)$ for some $v$ of
the form $v=u\cdot\delta^\prime$
with $u\in k^\prime$ as wanted.
\medskip

\noindent{\it Proof of the Main Theorem}. As we have already remarked,
we can suppose that the finitely generated group $\Gamma\subset\bar{\Q}^\times$
is the group of $S$-units in a number field $K$, where $K$ is Galois over
$\Q$ containing $\delta$  and  
 $S$ is stable under Galois conjugation as in lemmas 1,2,3.

Suppose by contradiction that we have
an infinite set $\Sigma\subset\Z\times\OS^\times$ of solutions $(q,u)$ to the inequality
(2.1). Let us define by induction a sequence $\{\delta_i\}\subset  K$,
an infinite  decreasing chain  $\Sigma_i$ of infinite subsets
of $\Sigma$ and an infinite strictly decreasing 
chain $k_i$ of subfields of $K$ with the following properties:
\smallskip

{\it For each natural number $n\geq 0$, 
$\Sigma_n\subset (\Z\times k_n)\cap\Sigma_{n-1}$,
$k_n\subset k_{n-1},\ k_n\neq k_{n-1}$, and
all but finitely many pairs $(q,u)\in\Sigma_n$  satisfy the inequalities
$|\delta_0\cdots\delta_n\cdot qu|>1$ and
$$
\|\delta_0\cdots\delta_n\cdot qu\|<q^{-d-\epsilon}H(u)^{-\epsilon/(n+1)}.\eqno(2.12)
$$
}
\smallskip

We shall eventually deduce a contradiction from the fact that the number field
$K$ does not admit any infinite decreasing chain of subfields.
We proceed as follows: put $\delta_0=\delta$,
$k_0=K\cap\R$ and $\Sigma_0=\Sigma$. Suppose we have defined
$\delta_n,k_n,\Sigma_n$ for a natural number $n$.  
Applying Lemma 3 with $k=k_n,\delta=\delta_0\cdots\delta_n$
we obtain that there exists an element $\delta_{n+1}\in k_n$, a proper subfield
$k_{n+1}$ of $k_n$ and an infinite set $\Sigma_{n+1}\subset\Sigma_n$
such that all pairs $(q,u)\in \Sigma_{n+1}$ verify $u=\delta_{n+1} v$ with $v\in k_{n+1}$.
Now, since for $v\in K$, $H(\delta_{n+1}v)\geq
H(\delta_{n+1})^{-1}H(v)$,  we have in particular that for  almost all  $v\in K$,
$H(\delta_{n+1}v)\geq H(v)^{(n+1)/(n+2)}$; then all but finitely many  such pairs satisfy   
$$
\|\delta_0\cdots\delta_n\delta_{n+1}\cdot qv\| <q^{-1-\epsilon}H(v)^{-\epsilon/(n+2)},
$$
and the inductive hypothesis is fullfilled.

The contradiction is then obtained as noticed above, concluding the proof of Main Theorem.
\bigskip

To prove Theorem 1 we also need the following 
\medskip

\noindent {\bf Lemma 4}. {\it Let $\alpha$ be an algebraic number. Suppose  that
for all $n$ in an infinite set $\Xi\subset\N$, there
exists a positive integer  $q_n\in\Z$ such that the sequence $\Xi\ni n\mapsto q_n$
 satisfies  
$$
\lim_{n\rightarrow\infty}{\log q_n\over n}=0\qquad {\rm and} \qquad 
{\rm Tr}_{\Q(\alpha)/\Q}(q_n\alpha^n)\in\Z\setminus\{0\}
$$
(the limit being taken for $n\in\Xi$). 
Then $\alpha$ is either the $h$-th root of a rational number
(for some  positive integer $h$) or  an algebraic integer.}
\medskip

\noindent{\it Proof}. It is essentially an application of Lemma
1, so it still depends on the Subspace Theorem. 
Let us suppose that $\alpha$ is not 
an algebraic integer. Let $K$ be the Galois
closure  of the extension $\Q(\alpha)/\Q$ and let $h$ be
the order of the torsion group $K^\times$. Since $\Xi$ is an infinite subset of
$\N$, there exists an integer $r\in\{0,\ldots,h-1\}$ such that  
infinitely many elements of $\Xi$  are of the form $n=r+hm$. Let us denote by
$\Xi^\prime$ the infinite subset of $\N$ composed of those integers $m$ such that
$r+hm\in \Xi$. Let $\sigma_1,\ldots,\sigma_d\in{\rm Gal}(K/\Q)$, where 
$d=[\Q(\alpha^h):\Q]$, be a set of
automorphisms of $K$ giving all the embedding $\Q(\alpha^h)\hookrightarrow K$.
In other terms, if $H\subset {\rm Gal}(K/\Q)$ is the subgroup fixing $\Q(\alpha^h)$,
the set $\{\sigma_1,\ldots,\sigma_d\}$ is a complete set of representatives of left
cosets of $H$ in ${\rm Gal}(K/\Q)$.
This proves that if $d=1$ then $\alpha$ is a $h$-th root of a rational number.
Suppose the contrary and suppose morevoer that
$\alpha$ is not an algebraic integer; we try to obtain a contradiction.
Since $\alpha$ is not an algebraic integer, there exists a finite absolute value $w$ of $K$
such that $|\alpha|_w>1$. Let $H^\prime\subset {\rm Gal}(K/\Q)$ be the subgroup fixing 
$\Q(\alpha)$, so that we have the  inclusions $H^\prime\subset H\subset {\rm Gal}(K/\Q)$,
corresponding to the chain $\Q(\alpha^h)\subset \Q(\alpha)\subset K$.

For each $i=1,\ldots,d$, let $T_i\subset{\rm Gal}(K/\Q)$ be a  complete
set of representatives for the set of automorphisms coinciding with
$\sigma_i$ in $\Q(\alpha^h)$, modulo $H^\prime$. In other words
the elements of $T_i$, when restricted to $\Q(\alpha)$, give all the
embedding of  $\Q(\alpha)\hookrightarrow K$ whose restrictions to $\Q(\alpha^h)$
coincide with $\sigma_i$. Also $T_1\cup\ldots\cup T_d$ is a complete set
of embeddings of $\Q(\alpha)$ in $K$.
Then we can write the trace of $\alpha^n$ ($n=r+hm$) as
$$
{\rm Tr}_{\Q(\alpha)/\Q}(\alpha^{r+hm})=
\sum_{i=1}^d\left(\sum_{\tau\in T_i}\tau(\alpha^r)\right)\sigma_i(\alpha^h)^m=
\lambda_1\sigma_1(\alpha^h)^m+\ldots+\lambda_d\sigma_d(\alpha^h)^m,
$$
where $\lambda_i=\sum_{\tau\in T_i}\tau(\alpha^r).$
Note that not all the coefficients $\lambda_i$ can vanish, since then
the trace would also vanish.
Since the trace of $q_n\alpha^n$ is integral, we have   
$|{\rm Tr}_{\Q(\alpha)/\Q}(\alpha^n)|_w\leq |q_n|_w^{-1}$;
then, since $\log q_n=o(n)$,  for every $\epsilon<\log|\alpha|_w/\log H(\alpha)$
and sufficiently large $m\in\Xi^\prime$ we have
$$
|\lambda_1\sigma_1(\alpha^{hm})+\ldots+\lambda_d\sigma_d(\alpha^{hm})|_w\leq
|q_{r+hm}|_w^{-1}<|\alpha^{hm}|_w\cdot H(\alpha^{hm})^{-\epsilon}.
$$
Applying Lemma 1 we arrive at a non-trivial equation of the form
$$
a_1\sigma_1(\alpha^h)^m+\ldots+a_d\sigma_d(\alpha^h)^m=0
$$
satisfied by infinitely many integers $m$.  
An application of the Skolem-Mahler-Lech theorem leads to the conclusion that 
for two indices $i\neq j$ 
some ratio  $\sigma_i(\alpha^h)/\sigma_j(\alpha^h)=(\sigma_i(\alpha)/\sigma_j(\alpha))^h$
 is a root of unity. 
But then it   equals $1$, since $\sigma_i(\alpha)/\sigma_j(\alpha)$ lies in $K^\times$ (and by
assumption $h$ is the exponent of the torsion group of $K^\times$). Then $\sigma_i$ coincides
with $\sigma_j$ on $\Q(\alpha^h)$. This contradiction concludes the proof. 
\medskip

\noindent{\it Proof of Theorem 1}. 
Let us suppose that the hypotheses of Theorem
1 are satisfied, so $\|\alpha\|^n<l^n$ for infinitely many $n$. Then, either
$\alpha$ is a $d$-th root of an integer, and we are done,
 or $\|\alpha^n\|\neq 0$ for all $n>1$. In this case, taking any 
$\epsilon <-\log l/\log H(\alpha)$ 
 we obtain that for infinitely many $n$
$$
0<\|\alpha^n\|<H(\alpha^n)^{-\epsilon}
$$
and of course the sequence $\alpha^n$ belongs to a finitely generated subgroup of
$\bar{\Q}^\times$. Then our Main Theorem, with $q=1, u=\alpha^n$, implies that
infinitely many numbers $\alpha^n$ are pseudo Pisot number, in particular have
non-zero integral trace. By Lemma 4, $\alpha^n$  is an algebraic integer,
so it is a Pisot number as wanted.

\bigskip

\noindent {\bf \S 3 Proof of Theorem 2.} 
We recall some basic facts about continued fractions.
Let $\alpha>1$ be a real irrational number. We will use the notation
$$
\alpha=[a_0,a_1,\ldots],
$$
where $a_0,a_1,\ldots$ are positive integers, to denote the continued fraction for 
$\alpha$. We also let $p_h,q_h$, for $h=0,1\ldots$ ,   be
the numerator and denominator of the truncated continued fraction
$[a_0,a_1,\ldots,a_h]$, so that, by a well known fact, for all $h$ we have
$$
\left|{p_h\over q_h}-\alpha\right|\leq {1\over q_h^2 \cdot a_{h+1}}.\eqno(4.1)
$$
Also, we have the recurrence relation, holding for $h\geq 2$,
$$
q_{h+1}=a_{h+1} q_h+q_{h-1}.\eqno(4.2)
$$
Let $T=\left(\matrix{a&b\cr c&d}\right)\in{\rm GL}_2(\Z)$ be a unimodular matrix 
and let $[b_0,b_1,\ldots]$ be the continued fraction of
$$
T(\alpha):={a\alpha+b\over c\alpha+d}.
$$
Then there exists an integer $k$ such that for all large $h$, $b_h=a_{h+k}$.
\smallskip

Consider now the case of a real {\it quadratic} irrational number $\alpha>0$, 
and let $\alpha^\prime$ be its  (algebraic) conjugate. 
The continued fraction of $\alpha$ is eventually periodic; it is purely periodic
if and only if $\alpha>1$ and $-1<\alpha^\prime<0$; we call such a quadratic
irrational {\it reduced}. In this case
the period of $-1/\alpha^\prime$ is the ``reversed" period of $\alpha$. 
For every quadratic irrational number there exists a unimodular transformation 
$T\in{\rm GL}_2(\Z)$ such that $T(\alpha)$ is reduced (this is equivalent to saying that
the expansion of $\alpha$ is eventually periodic).
Since the transformation $x\mapsto -1/x$ is also unimodular, it follows that
in any case $\alpha$ and $\alpha^\prime$ have a period of the same length.
We can summarize these facts as follows:
\medskip

\noindent {\bf Facts}. {\it Let $\alpha$ be a real quadratic irrational, 
$\alpha^\prime$ its conjugate.  
The continued fraction development of $\alpha$ is eventually periodic. 
The quadratic irrationals 
$$
|\alpha|,\ |\alpha^\prime|,\ \left|{1\over \alpha}\right|,\ 
\left|{1\over\alpha^\prime}\right|
$$
have periods of the same length.}
\medskip

We begin by proving  the easier part of Theorem 2, namely:
\smallskip

{\it If $\alpha>0$ is a unit in the quadratic ring $\Z[\alpha]$ then
the period for the continued fraction for $\alpha^n$ is uniformely bounded}.
\medskip

Let us begin with the case of odd powers
of a unit $\alpha>1$, with $\alpha^\prime<0$, i.e.  
a unit of norm $-1$. Let us denote by 
$t_n=\alpha^n+\alpha^{\prime n}={\rm Tr}_{\Q(\alpha)/\Q}(\alpha^n)$
the trace of $\alpha^n$. Then for all odd integers $n$, the number  $\alpha^n$
satisfies the relation
$(\alpha^n)^2-t_n\alpha^n-1=0$, i.e.
$$
\alpha^n=t_n+{1\over \alpha^n}
$$
hence its continued fraction is simply $[\overline{t_n}]$, and has period one.

We shall now consider units $\alpha^n$ of norm $1$, so including also even powers
of units of norm $-1$. Suppose $\alpha>1$, with $0<\alpha^\prime<1$ is
such a unit, so $\alpha^\prime=\alpha^{-1}$.
Denoting again by $t_n$ the trace of $\alpha^n$, we see that 
the integral part of $\alpha^n$ is $t_n-1$, so we put $a_0(n)=t_n-1$.
We have
$$
(\alpha^n-a_0(n))^{-1}=(1-\alpha^{\prime n})^{-1}={1-\alpha^n\over 
(1-\alpha^{\prime n})(1-\alpha^n)}=
{\alpha^n-1\over t_n-2}.
$$
For sufficiently large $n$, its integral part is $1$. So put $a_1(n)=1$
(at least for sufficiently large $n$). Then
$$
\left({\alpha^n-1\over t_n-2}-1\right)^{-1}=\left({\alpha^n-1-t_n+2\over t_n-2}\right)^{-1}
={t_n-2\over 1-\alpha^{\prime n}}.
$$
Again for sufficiently large $n$, the integral part of the above number is
$t_n-2$, so put $a_2(n)=t_n-2$.
Now
$$
{t_n-2\over 1-\alpha^{\prime n}}-(t_n-2)=
{t_n-2-t_n+t_n\alpha^{\prime n}+2-2\alpha^{\prime n}\over 1-\alpha^{\prime n}}=
{1+\alpha^{\prime 2n}-2\alpha^{\prime n}\over  1-\alpha^{\prime n}}=
{1-\alpha^{\prime n}}
$$
where we used the equation satisfied by $\alpha^{\prime n}$ to simplify the numerator.
Then
$$
\left({t_n-2\over 1-\alpha^{\prime n}}-a_2(n)\right)^{-1}= (1-\alpha^{\prime n})^{-1}=
{\alpha^n-1\over t_n-2}
$$
and the algorithm ends giving the continued fraction expansion, valid for sufficiently 
large $n$, 
$$
\alpha^n=[a_0(n),\overline{a_1(n),a_2(n)}]=[ t_n-1,\overline{1,t_n-2}] .
$$

\bigskip

We now prove  our last lemma:

\medskip

\noindent{\bf Lemma 5}. {\it  Let $\alpha>1$ be a real quadratic
number, not the square root of a rational number.
Let  $a_0(n),a_1(n),\ldots$ be the partial quotients of the
continued fraction for $\alpha^n$. 
Then either $\alpha$ is a Pisot number, or
for every $i$ with $i\neq 0$,
$$
\lim_{n\rightarrow \infty} {\log a_i(n)\over n}=0. \eqno(4.3)
$$
If $\alpha>1$ is the square root of a rational number, then (4.3) holds provided
the limit is taken over odd  integers $n$.
}
\medskip

{\it Proof}. We argue by contradiction. Let $h\in\N, h>0$ 
be the minimum index $i$ such that (4.3) does not hold. Then for 
a positive real $\epsilon$ and all $n$ in an infinite set  $\Xi\subset\N$
$$
a_h(n)>e^{\delta n}.\eqno(4.4)
$$
Here it is meant that $\Xi$ contains only odd integers if $\alpha$ is the
square root of a rational number.

On the other hand, since (4.3) holds for all $i=1,\ldots,h-1$, we have,
in view of the recurrence relation (4.2), that the denominators $q_{h-1}(n)$
of $p_{h-1}(n)/q_{h-1}(n)=[a_0(n),a_1(n),\ldots,a_{h-1}(n)]$ satisfy
$$
\lim_{n\rightarrow\infty} {\log q_h(n)\over n}=0.\eqno(4.5)
$$
Put $\displaystyle{ \epsilon={\delta\over 2\log H(\alpha)} }$. 
By the vanishing of the above limit we
have in particular
$$
q_{h-1}(n)^{1+\epsilon}\leq e^{{\delta n\over 2}}
$$
for sufficiently large $n\in \bar{\Xi}$. In view of (4.1) we have,
for all such $n$,
$$
\left|\alpha^n- {p_{h-1}(n)\over q_{h-1}(n)}\right|<
{1\over q_{h-1}(n)^2\cdot e^{\delta n}}
<{1\over q_{h-1}(n)^{3+\epsilon}\cdot e^{\delta n/2}}
$$
which can be rewritten as
$$
\|q_{h-1}(n)\alpha^n\|<q_{h-1}(n)^{-2-\epsilon}\cdot (H(\alpha^n))^{-\epsilon}.
$$
An application of our Main Theorem, with $u=\alpha^n$ (for $n\in \Xi$),
 gives the conclusion that 
for large $n$ the algebraic number $q_{h-1}(n)\alpha^n$ is pseudo-Pisot.
Now, taking into account (4.5),  
our Lemma 4 implies that   $\alpha$ is an algebraic integer
(it cannot have a vanishing trace, since otherwise it would
be the square root of a rational number, which we excluded).
This in turn implies that $\alpha^n$   is a Pisot number.
Since $\alpha$ is quadratic irrational and $\alpha^n$ is not a rational,
$\alpha$ itself is a Pisot number.

\medskip

We can now finish the proof of Theorem 2. 
Let $\alpha$ be a real quadratic 
irrational number; let us treat first the case that
$\alpha$ is not the root of a rational number. Suppose
 that for all positive integers $n$
in an infinite set $\Xi\subset\N$, the period of the continued fraction of
$\alpha^n$ has the same length $r$. We would like  to prove that  $\alpha$
is  a unit.

Let us first show that, by the above mentioned Facts, we can reduce to the case where   
$$
\alpha>1\qquad {\rm and}\qquad \alpha>|\alpha^\prime|.\eqno(4.6)
$$
In fact, after replacing if necessary $\alpha$ by $\pm 1/\alpha$ we can
suppose that $\alpha>1$.
Observe now that if   $\alpha=-\alpha^\prime$ then
$\alpha$ is the square root of a rational number, which we have excluded.
Then $|\alpha|\neq |\alpha^\prime|$ 
(since we cannot have $\alpha=\alpha^\prime$). If now $\alpha>|\alpha^\prime|$
we are done; otherwise 
$|\alpha^\prime|>\alpha>1$; in this case
we replace $\alpha$ by $\pm \alpha^\prime$ and obtain (4.6) as wanted.

So, from now on, let us suppose (4.6) holds for $\alpha$.
Then for all $n\geq n_0(\alpha)$ we have
$ \alpha^n>{\alpha^\prime}^n+2$
so there exists an integer $k_n$ such that
$$
\alpha^n-k_n>1, \quad -1<(\alpha^n-k_n)^\prime<0,
$$
so $\alpha^n-k_n$ is reduced. Then we have for all $n\in\Xi$,
$$
\alpha^n-k_n=[\overline{a_0(n),a_1(n),\ldots,a_{r-1}(n)}]
$$
where $a_0(n),\ldots,a_{r-1}(n)$ are positive integers and the 
period $(a_0(n),\ldots,a_{r-1}(n))$ is the same as the one for $\alpha^n$
(still for $n$ in the infinite set $\Xi$).
Suppose first that $\alpha$ is not a Pisot number. 
We would like to derive a contradiction. Our Lemma 5 
implies that all partial quotients  $a_i(n)$ satisfy (4.3), including
$a_0(n)=a_r(n)$. The algebraic numbers $\alpha^n-k_n$ and 
${\alpha^\prime}^n-k_n$ satisfy, for $n\in\Xi$, the algebraic equation
$$
x=[a_0(n);a_1(n),\ldots,a_{r-1}(n),x]={p_{r-1}(n)x+p_{r-2}(n)\over 
q_{r-1}(n)x+q_{r-2}(n)},
$$
which can be also written as
$$
 q_{r-1}(n) x^2+(q_{r-2}(n)-p_{r-1}(n))x-p_{r-2}(n)=0.
$$
 In view of (4.5), the coefficients of this equation have logarithms bounded by $o(n)$,
whence   the logarithmic heights of $x=\alpha^n-k_n$ and $x'={\alpha^\prime}^n-k_n$ are
$(o(n))$ for $n\in\Xi$. But then we would have
$$
\log|\alpha^n-{\alpha^\prime}^n|=o(n)
$$
which is clearly impossible.

We have then proved that $\alpha$ is a Pisot number, so it is
an algebraic integer and verifies $|\alpha^\prime|<1$. 
But now the quadratic irrational ${\pm 1\over \alpha^\prime}$ 
also satisfies (4.6), so the same reasoning implies that $1/\alpha^\prime$,
hence $1/\alpha$, is also an integer. This proves that $\alpha$ is a unit
as wanted.

The method is exactly the same in the  
  second case, when $\alpha=\sqrt{a/b}$ is the square root
of a positive rational number $a/b$. Here too, by replacing if necessary
$\alpha$ with $1/\alpha$, we reduce to the case $\alpha>1$.
In this case the pre-period has length one, so $(\alpha^n-k_n)^{-1}$ is 
reduced, where $k_n$ is the integral part of $\alpha^n$.
Under the hypothesis that the period length
of $\alpha^n$ remains bounded for an infinite set of odd integers,
we can still apply our Lemma 5 and conclude as before.\bigskip

\centerline{\title Appendix}\bigskip

In this section we show that  the word ``algebraic" cannot be removed from the
statement of Theorem 1. A little more precisely, we prove that: {\it There
exists a real number 
$\alpha >1$  such that $||\alpha^n||\le 2^{-n}$ for infinitely
many  positive integers $n$, but $\alpha^d$ is not a Pisot number, no matter the
positive integer $d$.}

Our assertion will follow from a simple construction which we are going to
explain. For a set $I$ of positive real numbers and a positive real $t$, we put
$I^t:=\{x^t: x\in I\}$.
\medskip

We start by choosing an arbitrary sequence $\{\beta_n\}_{n\ge 1}$  of real
numbers in the interval $[0,1/2]$; we shall then  define inductively integers
$b_0,b_1,\ldots $ and closed intervals $I_0,I_1,\ldots $ of real numbers, of
positive length and contained in $[2,\infty)$. We set 
$I_0=[2,3]$,
$b_0=1$
 and, having
constructed
$I_n, b_n$ we  continue as follows. Let
$b_{n+1}$ be a positive integer divisible by $n+1$ and such that $I_n^{b_{n+1}}$
contains an interval of length $>2$; such an integer will   certainly exist 
since by induction  every element of $I_n$ is $\ge 2$ and since $I_n$ has
positive length. 

Put $B_n=b_0\cdots b_n$. Then we choose  
$I_{n+1}$ to be any  subinterval of
$I_n^{b_{n+1}}$ of the shape 
$$
I_{n+1}=[q +\beta_{n+1}, q +\beta_{n+1}+2^{-B_{n+1}}],
$$
where $q=q_n$ is an integer. This will be possible since $I_n^{b_{n+1}}$ has
length $>2$. Clearly $q \ge 2$, so $I_{n+1}$ has the required properties.
\smallskip

For $n\in\N$, let now $J_n$ denote the closed interval $I_n^{1/B_n}$.
Since $I_{n+1}\subset I_n^{b_{n+1}}$, we have $J_{n+1}\subset J_n$ for all
$n$ and
$J_0=I_0=[2,3]$. Hence there exists $\alpha\in \bigcap_{n=0}^\infty J_n$ such
that $\alpha\ge 2$. 

Note that $\alpha^{B_n}\in I_n$; hence we conclude that: {\it For $n\ge 1$ the
fractional part
$\{\alpha^{B_n}\}$ of
$\alpha^{B_n}$ lies between $\beta_n$  and $\beta_n+2^{-B_n}$}. \medskip

To obtain our original assertion, we apply this construction by choosing 
$\beta_n=0$ for even integers $n\ge 1$ and $\beta_n=1/3$, say, for odd $n$. We
then have in the first place
$$
||\alpha^{B_n}||\le \{\alpha^{B_n}\}\le 2^{-B_n},
$$
for all even integers $n\ge 1$. This yields the first required property of
$\alpha$, because  $B_n\rightarrow \infty$. 

Further, suppose that $\alpha^d$ were a Pisot number, for some positive
integer
$d$. Then  we would have, for some fixed positive
$\lambda<1$ and integers $m\ge 1$,  
$||\alpha^{dm}||\ll \lambda^m$ (by a well-known
argument, used also in the previous proofs of this paper).
 However this is not possible, because for   odd integers $n>d$ we have that
$B_n$ is divisible by $d$ (recall $b_n$ is divisible by $n$ for $n>0$) and the
fractional part of $\alpha^{B_n}$ is between $1/3$ and $1/3+2^{-B_n}$, whence
the norm $||\alpha^{B_n}||$ is $\ge 1/6$.\medskip 

This concludes the argument. It will be noted that the construction yields more
precise conclusions; also, Theorem 1 obviously implies that any $\alpha$
likewise obtained is necessarily transcendental.

\bigskip

{\title References}\medskip

\item{[CZ]} - P. Corvaja, U. Zannier, On the length of the continued fraction 
             for the ratio of two power sums,  (2003)
            preprint NT/0401362 at  {\tt http://www.arXiv.org/} .\smallskip

\item{[M]} - K. Mahler, On the fractional parts of the powers of a rational
       number (ii), {\it Mathematika} {\bf 4} (1957), 122-124.\smallskip

\item{[MF]} - M. Mend\`es France, Remarks and problems on finite 
       and periodic continued fractions, {\it L'En\-seignement Math.} 
      {\bf 39} (1993)  249-257.\smallskip

\item{[S]} - W. M. Schmidt, Diophantine approximation and diophantine equations.
              {\it Lecture Notes in Mathematics } 1467,  Springer (1991).

\bigskip

\vfill

P. Corvaja\hfill U. Zannier

Dip. di Matematica e Informatica\hfill  Ist. Univ. Arch. - D.C.A.

Via delle Scienze\hfill S. Croce, 191

33100 - Udine (ITALY)\hfill  30135 Venezia (ITALY)

corvaja@dimi.uniud.it\hfill zannier@iuav.it

\end